# CORRECTION
# EFFICIENT PARAMETER ESTIMATION FOR SELF-SIMILAR PROCESSES

By Rainer Dahlhaus

*The Annals of Statistics* (1989) **17** 1749–1766

In the paper [1] the author claimed to have established asymptotic normality and efficiency of the Gaussian maximum likelihood estimator (MLE) for long range dependent processes (with the increments of the self-similar fractional Brownian motion as a special case—hence, the title of the paper). The case considered in the paper was Gaussian stationary sequences with spectral densities $f_\theta(x) \sim |x|^{-\alpha(\theta)} L_\theta(x)$, where $0 < \alpha(\theta) < 1$, $\theta \in \Theta \subset \mathbb{R}^p$ with $\Theta$ compact, and $L_\theta(x)$ being a slowly varying function at 0, a typical example being $f_\theta(x) = |x|^{-2d} h_\theta(x)$ and $\theta = (d, \theta_1, \ldots, \theta_{p-1})'$. The results were derived under a set of conditions termed (A0)–(A9), where in particular,

(A9) $\alpha$ is assumed to be continuous. Furthermore, there exists a constant $C$ with

$$|f_\theta(x) - f_{\theta'}(x)| \leq C |\theta - \theta'| f_{\theta'}(x) \tag{1}$$

uniformly for all $x$ and all $\theta, \theta'$ with $\alpha(\theta) \leq \alpha(\theta')$, where $|\cdot|$ denotes the Euclidean norm.

Unfortunately, this condition rules out long range dependent processes: Consider the case $f_\theta(x) = |x|^{-2d} h_\theta(x)$, $\theta = (d, \theta_1, \ldots, \theta_{p-1})'$, with $h_\theta(x)$ being bounded from above and below. Suppose $\theta'$ is fixed and $\theta = \theta_n$ with $\theta_n \to \theta'$, where $d = d_n = d' - 1/(2n)$. Then it is easy to verify that (A9) is not fulfilled for $x = 1/(n^n)$ and $n$ sufficiently large.

Luckily (A9) can be relaxed to

(A9$'$) $\alpha$ is assumed to be continuous.

We now briefly sketch how the results of the paper follow under this weaker condition. We assume that the reader is familiar with the notation and the results of the original paper.

REMARK. If $f_\theta(x) = |x|^{-\alpha(\theta)} h_\theta(x)$, the smoothness conditions in (A0)–(A8) on $f_\theta(x)$ are only fulfilled if $\alpha$ is three times differentiable [which is true, e.g., if $\alpha(\theta) = \theta_1$]. However, formally we need no such assumption since $\alpha$ appears only in upper bounds.

The most critical issue is the proof of consistency for the MLE and, in particular, the proof of the relation on page 1756, line 7, in [1]. This follows from equicontinuity in probability of the quadratic form

$$(2) \qquad Z_N^{(0)}(\theta) = \frac{1}{N}(\boldsymbol{X}_N - \mu_0 \boldsymbol{1})' T_N(f_\theta)^{-1} (\boldsymbol{X}_N - \mu_0 \boldsymbol{1}).$$

This is the quadratic form $Z_N^{(0)}(\theta)$ with $A_\theta^{(0)} = T_N(f_\theta)^{-1}$ in [1], Section 6.

PROPOSITION. *Suppose assumptions* (A0)–(A8) *in* [1] *and* (A9$'$) *hold. Then* $Z_N^{(0)}(\theta)$ *is equicontinuous in probability, that is, for each* $\eta > 0$ *and* $\varepsilon > 0$, *there exists a* $\delta > 0$ *such that*

$$(3) \qquad \limsup_{N \to \infty} \mathbf{P}\left( \sup_{|\theta_1 - \theta_2| \leq \delta} |Z_N^{(0)}(\theta_1) - Z_N^{(0)}(\theta_2)| > \eta \right) < \varepsilon.$$

PROOF. The proof follows exactly along the lines of the proof of Theorem 6.1 in [1]. The only additional argument needed (on page 1764, line 2) is the proof of $|\mathbf{E}S^k| \leq k!(2C)^k$ for $k = 1$, where $S = \{Z_N^{(0)}(\theta_1) - Z_N^{(0)}(\theta_2)\}_{jk}/|\theta_1 - \theta_2|$. This proof is quite technical. It can be found in the extended correction note [2]. □

With this result the following changes have to be made in [1]:

1. Page 1755, line 5 also holds without (A9);
2. Page 1756, line 7 follows from the equicontinuity stated above;
3. Lemma 5.5 can be skipped completely.

Furthermore, we mention that the proof of Lemma 6.2 contains a gap in that the proof of $|\mathbf{E}S^k| \leq k!(2C)^k$ is only given for $k > 1$. The case $k = 1$ is not straightforward. It can be obtained by a generalization of Lemma 2.1 of [2] to more functions. An alternative is to use [1], Theorem 5.1 [with the restriction $p(\beta - \alpha) < \frac{1}{2}$] and to note that (equi-) continuity of $Z_N^{(1)}(\theta)$ and $Z_N^{(2)}(\theta)$ is only needed at $\theta_0$ (in the proof of (iv) on page 1756 of [1]).

**Acknowledgments.** I am grateful to Professor Kung-Sik Chan for pointing out the error and for some discussions, and to Professor Hira Koul and two anonymous referees for several helpful comments.

INSTITUT FÜR ANGEWANDTE MATHEMATIK
IM NEUENHEIMER FELD 294
D-69120 HEIDELBERG
GERMANY
E-MAIL: dahlhaus@statlab.uni-heidelberg.de